\documentclass[review]{elsarticle}

\usepackage{lineno,hyperref}
\modulolinenumbers[5]

\usepackage{amsthm}
\usepackage{amssymb}
\usepackage{amsfonts}
\usepackage{amsmath}

\def\JJ{{\sf J}}
\def\Rr{{\sf R}}
\def\Hh{{\sf H}}

\def\RR{{\mathbb{R}}}
\def\NN{{\mathbb{N}}}
\def\ZZ{{\mathbb{Z}}}

\def\II{{\stackrel{\circ}{\mathrm{I}}}}

\def\f{{\sf {f}}}
\def\u{{\sf {u}}}
\def\v{{\sf {v}}}
\def\w{{\sf {w}}}
\def\h{{\sf {h}}}

\def\L{{\mathcal{L}}}
\def\C{{\mathcal{C}}}
\def\S{{\mathcal{S}}}

\def\I{{\mathrm I}}

\def\func#1#2{ \colon {#1} \longrightarrow {#2}}
\def\dim{{\rm dim\,}}

\newtheorem{thm}{Theorem}[section]
\newtheorem{lem}[thm]{Lemma} 
\newtheorem{propo}[thm]{Proposition} 
\newtheorem{corollary}[thm]{Corollary} 
\newtheorem{definition}[thm]{Definition}

\newdefinition{rmk}{Remark}
\newproof{pf}{Proof}
\newproof{pot}{Proof of Theorem \ref{thm2}}

\journal{arXiv}

\begin{document}

\begin{frontmatter}

\title{Explicit inverse of nonsingular Jacobi matrices}

\author{A.M. Encinas and M.J. Jim\'enez}
\address{Departament de Matem\`atiques, UPC, BarcelonaTech, Spain}

\begin{abstract}
We present here the necessary and sufficient conditions for the invertibility of tridiagonal matrices, commonly named Jacobi matrices, and explicitly compute their inverse. The techniques we use are related with the solution of Sturm--Liouville boundary value problems associated to second order linear difference equations. These boundary value problems can be expressed throughout a {\it discrete Schr\"odinger operator} and their solutions can be computed using recent advances in the study of linear difference equations. The conditions that ensure the uniqueness solution of the boundary value problem lead us to the invertibility conditions for the matrix, whereas the solutions of the boundary value problems provides the entries of the inverse matrix.
\end{abstract}

\begin{keyword}
tridiagonal matrices, second order linear difference equations, Sturm--Liouville boundary value problems, discrete Schr\"odinger operator, Chebyshev functions and polinomyals
\MSC[2010] 15B99\sep 31E05\sep  39A06
\end{keyword}

\end{frontmatter}

\linenumbers

\section{Preliminaries}

If we consider $n\in \mathbb{N}\setminus\{0\}$, the set $\mathcal{M}_n(\mathbb{R})$ of real matrices of size $n \times n$, and the sequences $a=\{a(k)\}_{k=0}^{n+1}\subset \mathbb{R}$, $b=\{b(k)\}_{k=0}^{n+1}\subset \mathbb{R}$, and $c=\{c(k)\}_{k=0}^{n+1}\subset \mathbb{R}$, then the Jacobi matrix associated with $a, b$ and $c$ is $\JJ(a,b,c)\in\mathcal{M}_{n+2}(\mathbb{R})$ given by
\begin{equation}
\label{matriz:Jirred}
\JJ(a,b,c)=\begin{bmatrix}
b(0) & -a(0) & 0&\cdots & 0 & 0\\
-c(0) & b(1) & -a(1) & \cdots & 0& 0\\
0 & -c(1) & b(2) & \cdots & 0 & 0\\
\vdots & \vdots &\vdots & \ddots & \vdots & \vdots \\
0 & 0 & 0 & \cdots & b(n) & -a(n)\\
0 & 0& 0& \cdots  & -c(n) & b(n+1)\end{bmatrix}
\end{equation}

Jacobi matrices appear frequently in both general Mathematics and Applied Mathematics, see \cite{M92}. As in this reference, we have chosen to write down the coefficients outside the main diagonal with negative sign. This is only a suitable convention, motivated by the existing relationship between Jacobi matrices and Schr\"odinger operators on a path, that we will use to analyze the invertibility of the Jacobi matrix. We must make also some assumptions about the coefficients of the matrix to avoid trivial situations or problems reducible to others with a minor order. Therefore, we will require $a(k)\not=0$ and $c(k)\not=0$, $k=0,\ldots,n$; since, in other case, $\JJ(a,b,c)$ is a reducible matrix and hence the inversion problem leads to the invertibility of a matrix of lower size. Moreover, the values of the coefficients for the sequences $a$ and $c$ at $n+1$ have no influence in the analysis of the matrix (\ref{matriz:Jirred}), since these coefficients do not appear in it. So, without loss of generality, we can impose $a(n+1)=c(n)$ and $c(n+1)=a(n)$. In the sequel, we also assume that $0^0=1$ and the usual convention that empty  sums and empty products are defined as $0$ and $1$, respectively.

The matrix $\JJ(a,b,c)$ is invertible if and only if for each $\f \in \mathbb{R}^{n+2}$ there exists $\u\in \mathbb{R}^{n+2}$ such that $\JJ(a,b,c)\u=\f$; which is equivalent to
\begin{equation}
\label{Jacobi:contorno}
\left\{\begin{array}{rll}
b(0)u(0)-a(0)u(1)=&\hspace{-.25cm}f(0),&\\[1ex]
-a(k)u(k+1)+b(k)u(k)-c(k-1)u(k-1)=&\hspace{-.25cm}f(k), &\hspace{-.5cm} k=1,\ldots,n,\\[1ex]
-c(n)u(n)+b(n+1)u(n+1)=&\hspace{-.25cm}f(n+1).&
\end{array}\right.
\end{equation}

Moreover, when this happens $\u$ is the unique solution of \eqref{Jacobi:contorno}. We can recognize in the previous equations the structure of a boundary value problem associated with a second order linear difference equation with coefficients $a,b,c$ and data $f$ or, equivalently, with a Schr\"odinger operator $\L_q$ with potential $q$ on the path $\I=\{0,\ldots,n+1\}$. Specifically, if $\II=\{1,\ldots, n\}$, $\delta(\I)=\{0,n+1\}$ and  $\C(\mathrm{I})$ is the vector space of real functions defined on $\mathrm{I}$, the Schr\"odinger operator with potential $q\in \C(\mathrm{I})$ on the path $\mathrm{I}$ is the linear operator $\L_q\func{\C(\mathrm{I})}{\C(\mathrm{I})}$ defined as
$$\left.\begin{array}{rll}
\L_q(u)(0)=&\hspace{-.25cm}a(0)\big(u(0)-u(1)\big)+q(0)u(0),&\\[1ex]
\L_q(u)(k)=&\hspace{-.25cm}a(k)\big(u(k)-u(k+1)\big)\\[1ex]
+&\hspace{-.25cm}  c(k-1)\big(u(k)-u(k-1)\big)+q(k)u(k),&\hspace{.25cm} k\in \II,\\[1ex]
\L_q(u)(n+1)=&\hspace{-.25cm}c(n)\big(u(n+1)-u(n)\big)+q(n+1)u(n+1),&
\end{array}\right\}$$
%
where $q\in \C(\mathrm{I})$ is defined as $q(0)=b(0)-a(0)$, $q(k)=b(k)-a(k)-c(k-1)$, $k\in \II$ and $q(n+1)=b(n+1)-c(n)$.
Identifying $\C(\mathrm{I})$ with $\RR^{n+2}$, and using this functional notation, Equation \eqref{Jacobi:contorno} is equivalent to the equation $\mathcal{L}_q(u)=f$ on $\mathrm{I}$; that is, to the 
 Sturm--Liouville value problem 
\begin{equation}
\label{Jacobi:contorno2}
\L_q(u)=f \hspace{.25cm}\hbox{on \,$\II$},\hspace{.25cm}\mathcal{L}_q(u)(0)=f(0)  \hspace{.25cm}\hbox{and}\hspace{.25cm}\mathcal{L}_q(u)(n+1)=f(n+1),
\end{equation}
where the identities $\L_q(u)=f$ on $\delta(\mathrm{I})$ play the role of boundary conditions, whereas the equation $\L_q(u)=f$ on $\II$  is named {\it the Schr\"odinger equation on $\II$ with data $f$}.

Therefore, the Jacobi matrix $\JJ(a,b,c)$ is invertible if and only if the linear operator $\L_q$ is invertible. In terms of the boundary value problem, the invertibility conditions for $\JJ(a,b,c)$ are exactly the same conditions to ensure that the boundary value problem is {\it regular}; that is, it has a unique solution for each given data and, hence, the computation of the inverse of $\JJ(a,b,c)$ can be reduced to the calculus of this solution. Implicitly or explicitly, determining the solutions for initial or final value problems for the Schr\"odinger equation is the strategy followed to achieve the inversion of tridiagonal matrices, see for instance \cite{FP01,FP05,SSH08,GK41,M01,MNNST98,U93,U94}; but either the general case is not analyzed, the explicit expressions of these solutions are not obtained, or the expressions obtained are excessively cumbersome.

\section{Initial value problems}

It is well-known that every {\it initial value problem} for the
Schr\"odinger equation on $\mathrm{\II}$ has a unique solution. Specifically, given  $f\in \C(\mathrm{I})$ and $m=0,\ldots,n$, for any $\alpha,\beta \in \RR$ there exists a unique $u\in \C(\mathrm{I})$ satisfying 
%
$$\L_q(u)=f\hspace{.25cm}\hbox{on $\II$}\hspace{.25cm}\hbox{and}\hspace{.25cm}u(m)=\alpha,\hspace{.25cm}u(m+1)=\beta.$$
%
In particular, when $m=n$, the above problem is also known as {\it final value problem}.

If $\S$ denotes the set of solutions of the homogeneous Schr\"odinger equation on $\II$ - that is $\L_q(u)=0$ on $\II$ - then $\S$ is a vector space such that $\dim \S=2$; while for any $f\in \C(\mathrm{I})$, the set $\S(f)$ of solutions of the 
Schr\"odinger equation on $\II$ with data $f$ satisfies $\S(f)\neq \emptyset$ and given $u \in \S(f)$, it is verified $\S(f)=u+\S$.

Given $u,v \in \C(\mathrm{I})$, their {\it Wronskian} or {\it Casoratian}, see \cite{A00}, is $w[u,v]\in \C(\mathrm{I})$ defined as
$$w[u,v](k)={\rm det}\begin{bmatrix}u(k) & v(k)\\u(k+1) & v(k+1)\end{bmatrix}=u(k)v(k+1)-v(k)u(k+1),\hspace{.25cm}0\leq k \leq n,$$
and as $w[u,v](n+1)=w[u,v](n)$. The Wronskian is a skew--symmetric bilinear form and and given $u,v\in \S$, either $w[u,v]=0$ or $w[u,v]\neq 0$ for any $k \in\II\cup\{0\} $. Moreover, $u$ and $v$ are linearly independent if and only if their Wronskian is non null and then $\{u,v\}$ form a basis of $\S$.

The {\it Green's
function} of the Schr\"odinger equation on $\II$ is the function \linebreak$g\in \mathcal{C}(\mathrm{I}\times \mathrm{I})$, defined  for any $s\in \mathrm{I}$
as $g(\cdot,s)$, the unique solution of the initial value problem
with conditions $g(s,s)=0$ and $g(s+1,s)=-\dfrac{1}{a(s)}$, when $0\le s\le n$,  and as the unique solution of the initial value problem
with conditions $g(n+1,n+1)=0$ and $g(n,n+1)=\dfrac{1}{a(n+1)}$ when $s=n+1$. Notice that $g(s,s+1)=\dfrac{1}{c(s)}$ for any $s=0,\ldots,n$. Therefore, if for any $s=0,\ldots,n$ we consider $u=g(\cdot,s)$ and $v=g(\cdot,s+1)$, then 
$$w[u,v](s)=g(s,s)g(s+1,s+1)-g(s+1,s)g(s,s+1)=\dfrac{1}{a(s)c(s)}$$
which implies that $\{g(\cdot,s),g(\cdot,s+1)\}$ is a basis of $\S$.
Moreover, for any $f\in \C(\mathrm{I})$ and $m=0,\ldots,n$, the function $u\in \C(\mathrm{I})$ given by 
%
$$u(k)=\sum\limits_{s=\min\{k,m\}+1}^{\max\{k,m\}}g(k,s)f(s),\hspace{.25cm}k\in \mathrm{I}$$
%
is the unique solution of the initial value problem $\L_q(u)=f$ on $\II$, and $u(m)=u(m+1)=0$.

It will be very useful to introduce the {\it companion function} defined as
$$\rho(k)=\prod\limits_{s=0}^{k-1}\dfrac{a(s)}{c(s)},\hspace{.5cm}k=0,\ldots,n+1.$$ 
Notice that $\rho(0)=1$.

Remembering the assumption $a(k),c(k)\not=0$, $0\leq k \leq n$, it is easy to prove that $\rho(k)a(k)=\rho(k+1)c(k)$. Moreover, the companion function verifies the following meaningful result.

\begin{propo}\label{cte}
Given $u,v \in \C(\mathrm{I})$, then 
$$a(k)w[u,v](k)=c(k-1)w[u,v](k-1) \hspace{.25cm}\hbox{for any}\hspace{.25cm} k \in \,\stackrel{\circ}{\mathrm{I}}.$$
Therefore, the multiplication of functions $\rho a w[u,v]$ is constant in $\mathrm{I}$ and is zero if and only if $u$ and $v$ are linearly dependent.
\end{propo}

\section{Regular Sturm--Liouville boundary value problems}

A boundary condition at $0$ is a linear function $\mathfrak{c}\func{\C(\mathrm{I})}{\RR}$ of the form $\mathfrak{c}(u)=\alpha u(0)+\beta u(1)+\gamma u(n)+\delta u(n+1)$, and a boundary condition at $n+1$ is a linear function $ \mathfrak{d}\func{\C(\mathrm{I})}{\RR}$ of the form $\mathfrak{d}(u)=\hat \alpha u(0)+\hat\beta u(1)+\hat \gamma u(n)+\hat \delta u(n+1)$. The pair $(\mathfrak{c},\mathfrak{d})$ is named Sturm--Liouville conditions if $\gamma=\delta=\hat \alpha=\hat \beta=0$, see \cite{BCE09}. Therefore, defining the pair of Sturm--Liouville conditions $(\mathfrak{c}_1,\mathfrak{c}_2)$ as 
$$\begin{array}{ll}
\mathfrak{c}_1(u)=&\hspace{-.25cm}\mathcal{L}_q(u)(0)=b(0)u(0)-a(0)u(1),\\
 \mathfrak{c}_2(u)=&\hspace{-.25cm}\mathcal{L}_q(u)(n+1)=-c(n)u(n)+b(n+1)u(n+1),\end{array}$$ 
and according to Equation (\ref{Jacobi:contorno2}), we must consider the Sturm--Liouville boundary value problem $(\L_q,\mathfrak{c}_1,\mathfrak{c}_2)$; that is, for any $f\in {\mathcal C}(\I)$,
we should determine if there exists $u\in {\mathcal C}(\I)$ such that
$${\mathcal
L}_q(u)=f\hspace{.25cm}\hbox{on}\hspace{.25cm}\stackrel{\circ}{\I},\hspace{.5cm}\mathfrak{c}_1(u)=f(0)\hspace{.25cm}\hbox{and} \hspace{.25cm}\mathfrak{c}_2(u)=f(n+1).$$
The boundary value problem $(\L_q,\mathfrak{c}_1,\mathfrak{c}_2)$ is called homogeneous when $f=0$.

We are only interested in {\it regular problems}; that is, in those boundary value problems with a unique solution. For the resolution of this sort of boundary value problems, we use the so--called {\it resolvent kernel}, see \cite[Sections 2 and 3]{EnJi17}, and the process of determining the resolvent kernel always depends on an appropriate choice of solutions of the corresponding homogeneous Schr\"odinger equation. We reproduce here some of the main results of the above--mentioned work of the authors, essential for the main result developed in the next section. Therefore, for more details or to check out proofs that are not included on the present section, see \cite{EnJi17}.

If $g$ is the Green function of the Schr\"odinger equation on $\II$, the value 
$$D_{a,b,c}=
\mathfrak{c}_1\big(g(\cdot,0)\big)\mathfrak{c}_2\big(g(\cdot,1)\big)-\mathfrak{c}_2\big(g(\cdot,0)\big)\mathfrak{c}_1\big(g(\cdot,1)\big)$$
encompasses information of both the Schr\"odinger equation on $\II$ and the pair of boundary conditions $(\mathfrak{c}_1,\mathfrak{c}_2)$. In fact, we next show that  it plays a fundamental  role in the analysis of the Sturm--Liouville problem.

\begin{definition} The boundary value problem $(\L_q,\mathfrak{c}_1,\mathfrak{c}_2)$ is called {\rm regular} if the solution of the corresponding homogeneous problem is unique, and so the null one.\end{definition}

\begin{propo} The following assertions are equivalent:
\begin{itemize}
\item[(i)]  The boundary value problem $(\L_q,\mathfrak{c}_1,\mathfrak{c}_2)$ is regular.
\item[(ii)] For any $f\in {\mathcal C}(\I)$ the corresponding boundary value problem has a solution (and hence a unique solution).
\item[(iii)] $D_{a,b,c}\not=0$.
\end{itemize}
\end{propo} 
\begin{pf}  If $z_1=g(\cdot,0)$ and $z_2=g(\cdot,1)$, then $\{z_1,z_2\}$ form a  basis of solutions of the homogeneous Schr\" odinger equation ${\mathcal L}_q(u)=0$ on $\II$. If given $f\in \mathcal{C}(\I)$ we consider  
a solution $y$ of the Schr\"odinger equation 
with data $f$ on $\II$, the expression $u=\alpha z_1+\beta z_2+y$ where $\alpha,\beta\in \RR$, determines all the solutions of the Schr\"odinger equation on $\II$. Therefore, $u=\alpha z_1+\beta z_2+y$ denotes a solution of the boundary value problem
$${\mathcal
L}_q(u)=f\hspace{.25cm}\hbox{on}\hspace{.25cm}\stackrel{\circ}{\I},\hspace{.5cm}\mathfrak{c}_1(u)=f(0)\hspace{.25cm}\hbox{and} \hspace{.25cm}\mathfrak{c}_2(u)=f(n+1),$$
if and only if $\alpha$ and $\beta$ are solutions of the linear system
$$\left[\begin{array}{cc}\mathfrak{c}_1(z_1)& \mathfrak{c}_1(z_2)\\\mathfrak{c}_2(z_1)&\mathfrak{c}_2(z_2)\end{array}\right]
\left[\begin{array}{c}\alpha\\\beta\end{array}\right]=\left[\begin{array}{c}f(0)-\mathfrak{c}_1(y)\\f(n+1)-\mathfrak{c}_2(y)\end{array}\right].
$$ 
When $f$ goes over $\mathcal{C}(\I)$, then the right term of the previous system goes over the whole $\RR^2$. Therefore, the system has a solution for any $f\in \C(\I)$ if and only if the coefficient matrix is non--singular and, hence, the system has a unique solution. As the homogeneous system associated with the previous one determines the solutions of the homogeneous boundary value problem, the problem is regular if the homogeneous system has as its unique solution the null one. Therefore, (i) and (ii) are equivalent and, in addition, the coefficient matrix is non--singular and it implies that its determinant is different from $0$. Hence, (i) and (iii) are equivalent.\qed
\end{pf}

In the sequel, for any $s\in \I$, we denote by $\varepsilon_s\in \C(\I)$ the {\it Dirac function at $s$}. Therefore $\varepsilon_s(s)=1$ and $\varepsilon_s(k)=0$, when $k\not=s$.

\begin{definition}\label{green:contorno} Let $(\L_q,\mathfrak{c}_1,\mathfrak{c}_2)$ be a regular boundary value problem. We call {\rm resolvent kernel of the boundary value problem} to $R_{a,b,c}\func{\I\times \I}{\RR}$ characterized by
$$
\L_q\big(R_{a,b,c}(\cdot,s)\big)=\varepsilon_s\hspace{.15cm}\hbox{on}\hspace{.15cm}\stackrel{\circ}{\I}, \hspace{.15cm}
\mathfrak{c}_1\big(R_{a,b,c}(\cdot,s)\big)=\varepsilon_{s}(0),\hspace{.15cm}\mathfrak{c}_2\big(R_{a,b,c}(\cdot,s)\big)=\varepsilon_{s}(n+1)
$$
for any $s\in \I$.
\end{definition}

Notice that for any $s\in \I$, $R_{a,b,c}(\cdot,s)$ is the unique solution of the Sturm-Liouville problem for the data $\varepsilon_s$ and hence it makes sense when the boundary value problem is regular. The role of the resolvent kernel is showed in the following result.

\begin{propo}
\label{resolvente:solucion} If the boundary value problem $(\L_q,\mathfrak{c}_1,\mathfrak{c}_2)$ is regular and $R_{a,b,c}$  is the resolvent kernel, then for any $f\in {\mathcal C}(\I)$ the function

$$u(k)=\sum_{s \in\I} R_{a,b,c}(k,s)\,f(s),\hspace{.25cm}k\in \I,$$
is the unique solution of the boundary value problem with data $f$, i.e.
$$\L_q(u)=f\hspace{.25cm}\hbox{on}\hspace{.25cm}\stackrel{\circ}{\I}, \hspace{.25cm}\mathfrak{c}_1(u)=f(0),\hspace{.25cm}\mathfrak{c}_2(u)=f(n+1).$$
\end{propo}

\begin{definition}
\label{fundamentales}
We call {\rm fundamental solutions of the homogeneous Schr\"odinger equation on $\stackrel{\circ}{\I}$, related to the boundary conditions $\mathfrak{c}_1$ and $\mathfrak{c}_2$} or, simply, {\rm fundamental solutions}, to $\Phi_{a,b,c},\Psi_{a,b,c}\in \mathcal{C}(\I)$, the unique solutions of the homogeneous Schr\"odinger equation on $\stackrel{\circ}{\I}$ determined respectively by the conditions
$$\begin{array}{rlrl}
\Phi_{a,b,c}(0)=&\hspace{-.25cm}a(0),& \hspace{.25cm}\Phi_{a,b,c}(1)=&\hspace{-.25cm}b(0),\\[1ex]
\Psi_{a,b,c}(n)=&\hspace{-.25cm}b(n+1),&\hspace{.25cm}\Psi_{a,b,c}(n+1)=&\hspace{-.25cm}c(n).
\end{array}$$
\end{definition}

Notice that $\Phi_{a,b,c}$ is the solution of a initial value problem, whereas $\Psi_{a,b,c}$ is the solution of a final value problem. The reason to choose these definitions for the fundamental solutions is shown in the following result. 

\begin{propo}
\label{contorno:base}
If $\Phi_{a,b,c}$ and $\Psi_{a,b,c}$ are the fundamental solutions of the homogeneous Schr\"odinger equation on $\stackrel{\circ}{\I}$, related to the boundary conditions $\mathfrak{c}_1$ and $\mathfrak{c}_2$, then 
$\mathfrak{c}_1(\Phi_{a,b,c})=\mathfrak{c}_2(\Psi_{a,b,c})=0$, $\mathfrak{c}_2(\Phi_{a,b,c})=a(0)c(0)D_{a,b,c}$. 
Moreover, 
$$\mathfrak{c}_1(\Psi_{a,b,c})=c(0)a(n)\rho(n)D_{a,b,c}=w[\Psi_{a,b,c},\Phi_{a,b,c}](0).$$
\end{propo}

\begin{pf} 
Consider $\{u,v\}$ the basis of solutions of the homogeneous Schr\"odinger equation satisfying $u(0)=1$, $u(1)=0$, $v(0)=0$ and $v(1)=1$; that is, \linebreak$u=c(0)g(\cdot,1)$ and $v=-a(0)g(\cdot,0)$. 
Moreover, $w[u,v](0)=1$.

If we prove that
$$\Phi_{a,b,c}=\mathfrak{c}_1(u)v-
\mathfrak{c}_1(v)u\hspace{.5cm}
\hbox{and}\hspace{.5cm}\Psi_{a,b,c}=
a(0)^{-1}a(n)\rho(n)\big(\mathfrak{c}_2(v)u -
\mathfrak{c}_2(u)v\big),$$
then, clearly, $\mathfrak{c}_1(\Phi_{a,b,c})=
\mathfrak{c}_2(\Psi_{a,b,c})=0$,
$\mathfrak{c}_2(\Phi_{a,b,c})=a(0)c(0)D_{a,b,c}$ and
$$\mathfrak{c}_1(\Psi_{a,b,c})=
a(0)^{-1}a(n)\rho(n)\mathfrak{c}_2(\Phi_{a,b,c})=c(0) a(n)\rho(n)D_{a,b,c}.$$
Moreover,
$$\begin{array}{rl}
w[\Phi_{a,b,c},\Psi_{a,b,c}](0)=&\hspace{-.25cm}a(0)^{-1}a(n)
\rho(n)\big(\mathfrak{c}_1(v)
\mathfrak{c}_2(u)-\mathfrak{c}_1(u)
\mathfrak{c}_2(v)\big)\\[1ex]
=&\hspace{-.25cm} -c(0)a(n)\rho(n)D_{a,b,c}.\end{array}$$

To end the proof, let us consider the functions 
$$z=\mathfrak{c}_1(u)v-\mathfrak{c}_1(v)u\hspace{.5cm}\hbox{and}\hspace{.5cm}\hat z=a(0)^{-1}a(n)\rho(n)\big(
\mathfrak{c}_2(v)u-\mathfrak{c}_2(u)v\big).$$

Then  $z(0)=-\mathfrak{c}_1(v)=a(0)$, $z(1)=\mathfrak{c}_1(u)=b(0)$ and on the other hand,
$$\begin{array}{rl}
\hat z(n)=&\hspace{-.25cm}a(0)^{-1}a(n)\rho(n)\big(
\mathfrak{c}_2(v)u(n)-\mathfrak{c}_2(u)v(n)\big)\\[1ex]
=&\hspace{-.25cm}b(n+1)a(0)^{-1}a(n)\rho(n)w[u,v](n)\\[1ex]
=&\hspace{-.25cm}b(n+1)a(0)^{-1}a(0)\rho(0)w[u,v](0)=b(n+1),\\[1ex]
\hat z(n+1)=&\hspace{-.25cm}a(0)^{-1}a(n)\rho(n)\big(
\mathfrak{c}_2(v)u(n+1)-\mathfrak{c}_2(u)v(n+1)\big)\\[1ex]
=&\hspace{-.25cm}c(n)a(0)^{-1}a(n)\rho(n)w[u,v](n)\\[1ex]
=&\hspace{-.25cm}c(n)a(0)^{-1}a(0)\rho(0)w[u,v](0)=c(n).\end{array}$$
The uniqueness of the solution of any initial value problem concludes that \linebreak$z=\Phi_{a,b,c}$ and $\hat z=\Psi_{a,b,c}$.\qed
\end{pf}

\begin{corollary}
The boundary value problem $(\L_q,\mathfrak{c}_1,\mathfrak{c}_2)$ is regular if and only if the fundamental solutions are a basis of solutions of the homogeneous Schr\"odinger equation on $\stackrel{\circ}{\I}$. 
\end{corollary} 

The next step in this section is to obtain  the resolvent kernel for a regular boundary value problem with Sturm-Liouville conditions in terms of the fundamental solutions, see \cite{EnJi17} for its proof.

\begin{thm}
\label{resolvente}
The Sturm--Liouville boundary value problem $(\L_q,\mathfrak{c}_1,\mathfrak{c}_2)$ is re\-gu\-lar if and only if $b(0)\Psi_{a,b,c}(0)\neq a(0)\Psi_{a,b,c}(1)$ or, equivalently, iff $c(n)\Phi_{a,b,c}(n)\neq \linebreak b(n+1)\Phi_{a,b,c}(n+1)$ and its resolvent kernel is determined by 
$$\begin{array}{rl}
R_{a,b,c}(k,s)=&\hspace{-.25cm}\displaystyle  \dfrac{\Phi_{a,b,c}(\min\{k,s\})\Psi_{a,b,c}(\max\{k,s\})}{a(0)\big[b(0)\Psi_{a,b,c}(0)-a(0)\Psi_{a,b,c}(1)\big]}\rho(s),
\end{array}$$
for any $k,s=0,\ldots,n+1$.
\end{thm}

Finally, let us remind that the boundary conditions associated with the Jacobi matrix were $\mathfrak{c}_1(u)=\mathcal{L}_q(u)(0)$ and $\mathfrak{c}_2(u)=\mathcal{L}_q(u)(n+1)$, so the boundary value problem $(\L_q,\mathfrak{c}_1,\mathfrak{c}_2)$ associated with the inversion of that matrix is the Poisson equation  $\L_q(u)=f$ on $\I$. Applying now Theorem \ref{resolvente} to this equation, we obtain the fundamental result for the inversion of Jacobi matrices. 

\begin{corollary}
\label{Poisson:resolvente} 
The Schr\"odinger operator $\L_q$ is invertible if and only if \linebreak $b(0)\Psi_{a,b,c}(0)\not=a(0)\Psi_{a,b,c}(1)$ and, moreover, given $f\in \C(\I)$, 
$$ 
(\L_q)^{-1}(f)(k)= \sum_{s \in\I}\dfrac{\Phi_{a,b,c}(\min\{k,s\})\Psi_{a,b,c}(\max\{k,s\})}{a(0)\Big[b(0)\Psi_{a,b,c}(0)-a(0)\Psi_{a,b,c}(1)\Big]}\,\rho(s) f(s),$$
for any $k=0,\ldots,n+1$.
\end{corollary}

\section{The inverse of a Jacobi matrix}

The invertibility conditions of the Jacobi matrix  $\JJ(a,b,c)$ described in Equation (\ref{matriz:Jirred}), as well as determining its inverse $\JJ^{-1}=\Rr=(r_{ij})$ in terms of the  solutions $\Phi_{a,b,c}$ and $\Psi_{a,b,c}$ of the Schr\"odinger equation, are described in Corollary \ref{Poisson:resolvente}. So, to obtain the explicit values of the entries of $\Rr$, the next step is to compute explicitly the functions $\Phi_{a,b,c}$ and $\Psi_{a,b,c}$, that can be seen as the solutions of an initial and a final value problem respectively, associated with the second order linear difference equation with coefficients $a,b$ and $c$ that corresponds to the Schr\"odinger equation. To compute these solutions we will use recent advances in the study of difference equations developed by the authors in \cite{EJ17}. In particular, in Section 7 of this work it has been proved that the solution of any initial value problem for a second order difference equation with any data $f$, can be expressed as a linear combination of the functions $P_k(x,y)$ called {\it $k$-th Chebyshev functions} and defined for any $x,y\in \C(\ZZ)$ as
\begin{equation}
\label{Ch}
P_0(x,y)=1,\hspace{.15cm} P_{-1}(x,y)=0\hspace{.15cm}\hbox{and}
\hspace{.15cm} P_k(x,y)=\sum\limits_{m=0}^{\lfloor \frac{k}{2}\rfloor}(-1)^m\sum\limits_{\alpha\in \ell_{k}^m}x^{\bar \alpha}y^\alpha,\hspace{.15cm}k\geq 1.
\end{equation}
We reproduce here some brief explanations about the notation involved in Equation (\ref{Ch}), for the sake of completeness. 
The parameter $\alpha=(\alpha_1,\ldots,\alpha_p)$ is a {\it  binary multi--index of order p}; i.e. $\alpha$ is a $p$--tuple $\alpha=(\alpha_1,\ldots,\alpha_{p})\in \{0,1\}^p$, and its length is defined as $|\alpha|=\sum\limits_{j=1}^{p}\alpha_j\le p$. Given $\alpha\in \{0,1\}^p$ and a function $a\in \C(\ZZ)$, we consider the value $a^\alpha=\prod\limits_{j=1}^{p}a(j)^{\alpha_j}$.  Given $p\in \NN\setminus\{0\}$, we denote by $i_1,\ldots,i_m$ the indices such that $1\le i_1<\cdots<i_m\le p$ and $\alpha_{i_j}=1$, $j=1,\ldots,m$. We just need to consider the binary multi--indexes $\alpha$ of order $p$ in the set $\ell_p$ defined as
\begin{enumerate}[(i)]{\parskip=.2cm
\item $\ell_p^0=\{\alpha: |\alpha|=0\}=\{(0,\ldots,0)\}$, for  $p\in \NN\setminus\{0\}$, 
\item $\ell_p^1=\{\alpha: \alpha_p=0\hspace{.15cm}\hbox{and}\hspace{.15cm}|\alpha|=1\}$, for $p\ge 2$, 
\item  
$\ell_p^m=\{\alpha:\alpha_p=0,\hspace{.15cm} |\alpha|=m\hspace{.15cm}\hbox{and}\hspace{.15cm}i_{j+1}-i_j\ge 2,\hspace{.15cm}j=1,\ldots,m-1\}$, for $p\ge 4$ and $m=2,\ldots,\lfloor\frac{p}{2}\rfloor$.
}
\end{enumerate}

Finally, $\bar \alpha$ is the binary multi--index of the same order as $\alpha$ defined by
$$\bar \alpha_{i_j}=\bar \alpha_{i_j+1}=0, \hspace{.15cm}j=1,\ldots,m,\hspace{.25cm}\hbox{and}\hspace{.25cm}\bar \alpha_i=1\hspace{.25cm}\hbox{otherwise}.$$

The name of Chebyshev function for (\ref{Ch}) is justified due to its relation with the usual Chebyshev polynomials of second kind, since $P_k(x,y)$ can be identified with them when $x$ and $y$ are constant sequences. In that case, $P_0(x,y)=1$ and $P_{-1}(x,y)=0$ and moreover, since $\#\ell_k^m={k-m\choose m}$  for any $k\in \NN^*$, we obtain that 
%
$$P_k(x,y)=\sum\limits_{m=0}^{\lfloor \frac{k}{2}\rfloor}(-1)^m{k-m\choose m}x^{k-2m}y^m.$$
%
Clearly, for any $k\ge -1$ and any constant sequence $x$, we have 
%
$$U_k(x)=P_k(2x,1)=\sum\limits_{m=0}^{\lfloor \frac{k}{2}\rfloor}(-1)^m{k-m\choose m}(2x)^{k-2m},$$
%
that is known as the standard {\it $k$--{\it th} Chebyshev polynomial of second kind}, see \cite{ABD05} and also \cite{BCE09,EJ16}. Definitely, for constant sequences $x$ and $y$, it is satisfied
$$P_k(x,y)=y^{\frac{k}{2}} \sum\limits_{m=0}^{\lfloor \frac{k}{2}\rfloor}(-1)^m{ k-m \choose m}\left(\dfrac{x}{\sqrt{y}}\right)^{k-2m}=y^{\frac{k}{2}}U_k\left(\dfrac{x}{2\sqrt{y}}\right),\hspace{.25cm}k\ge 1.$$

Now we are ready to compute the basis of solutions $\{\Phi_{a,b,c}(k),\Psi_{a,b,c}(k)\}$ of the homogeneous Schr\"odinger equation $\L_{q}(u)=0$ on $\I$ applying the results showed in \cite{EJ17} on second order difference equations, that is through a linear combination of the Chebyshev functions $P_k(b,ac)$ and $P_k(b_m,a_mc_m)$,  where $a,b,c \in \C(\I)$ are the coefficients of the second order difference equation associated to the Schr\"odinger equation, and given $a\in\C(\ZZ)$ and $m \in \NN$, the function $a_{m}$ corresponds to the $m$--shift of $a$, so $a_m=a(k+m)$. We must consider, for this first result and most of those that will appear from now on, the functions $\Phi_{_\JJ},\Psi_{_\JJ}\in \C(\I)$ defined as
$$\begin{array}{rlr}
\Phi_{_\JJ}  (0)=&\hspace{-.25cm}1, & \\[1ex]
\Phi_{_\JJ}  (k)=&\hspace{-.25cm} b(0)P_{k-1}(b,ac)-a(0)c(0)P_{k-2}(b_1,a_1c_1),&\hspace{-.65cm}  k=1\ldots,n+1,\\[2ex]
\Psi_{_\JJ} (k)=&\hspace{-.25cm} b(n+1)P_{n-k}(b_k,a_kc_k)-a(n)c(n)P_{n-k-1}(b_k,a_kc_k)     , &\hspace{-.65cm}  k=0,\ldots,n,\\[2ex]
\Psi_{_\JJ} (n+1)=&\hspace{-.25cm}1, &
\end{array} $$
and the value 
$$\begin{array}{rl}
D_{_\JJ}=&\hspace{-.25cm}b(0)\Big[b(n+1)P_{n}(b,ac)   -a(n)c(n)P_{n-1}(b,ac)\Big]\\[1ex]
-&\hspace{-.25cm}a(0) c(0)
\Big[b(n+1)P_{n-1}(b_1,a_1c_1)-a(n)c(n)P_{n-2}(b_1,a_1c_1)   \Big]. 
\end{array}$$

\begin{lem}
\label{phi_jacobi}
For any $k=0,\ldots,n+1$, it is satisfied that
$$\Phi_{a,b,c}(k)=a(0)\Big(\prod\limits_{s=0}^{k-1}a(s)\Big)^{-1}\Phi_{_\JJ}(k)\hspace{.5cm}\hbox{and}\hspace{.5cm} \Psi_{a,b,c}(k)=c(n)\Big(\prod\limits_{s=k}^{n}c(s)\Big)^{-1}\Psi_{_\JJ}(k)$$
and, moreover, 
$$b(0)\Psi_{a,b,c}(0)-a(0)\Psi_{a,b,c}(1)=D_{_\JJ}\Big(\prod\limits_{s=0}^{n-1}c(s)\Big)^{-1}.
$$
\end{lem}

\begin{pf}
Applying \cite[Theorems 4.3 and 7.4]{EJ17}, we have that $\Phi_{a,b,c}$ is a linear combination of the Chebyshev functions $\{P_{k-2}(b,ac),P_{k-1}(b_1,a_1c_1)\}$ and, in addition, $\Psi_{a,b,c}$ is a linear combination of $\{P_{n-k-1}(b_k,a_kc_k),P_{n-k}(b_k,a_kc_k)\}$.
To obtain all the results, we must just to impose the conditions
$$\Phi_{a,b,c}(0)=a(0),\hspace{.15cm}\Phi_{a,b,c}(1)=b(0),\hspace{.15cm}\Psi_{a,b,c}(n)=b(n+1),\hspace{.15cm}\Psi_{a,b,c}(n+1)=c(n).$$
\qed
\end{pf}

Before showing the main result for the explicit inversion of a Jacobi matrix, we add a previous result extracted from \cite[Theorem 3.3]{MNNST98}, that allow us to compute also the determinant of the inverse matrix.

\begin{lem} 
\label{Green:generalizada}
If $\Rr=(r_{ij})\in \mathcal{M}_m(\RR)$ is an irreducible and invertible matrix, the following statements are equivalents:
\begin{itemize}
\item[(i)] There exists a diagonal and invertible matrix $\Hh=(h_{j})$ such that $\Rr\Hh^{-1}$ is a Green's matrix; that is, there exist $\v,\w, \h\in \RR^m$, where $h_j\not=0$,\linebreak $j=1,\ldots,m$, such that 
$$r_{ij}=h_jv_{\min\{i,j\}}w_{\max\{i,j\}}=\left\{\begin{array}{cl} v_iw_jh_j, & \hbox{si $i\le j$},\\[1ex]
v_jh_jw_i; & \hbox{si $i\ge j$},\end{array}\right.$$
that is, 
$$\Rr=\begin{bmatrix}v_1 & v_1& v_1&\cdots & v_1\\
h_1v_1 & v_2 & v_2& \cdots & v_2\\
h_1v_1 & h_2v_2 & v_3 & \cdots & v_2\\
\vdots & \vdots& \vdots &\ddots & \vdots\\
h_1v_1 & h_2v_2 & h_3v_3&  \cdots &v_m \end{bmatrix}\circ \begin{bmatrix}h_1w_1 & h_2w_2& h_3w_3 & \cdots & h_nw_n\\
w_2 & h_2w_2 & h_3w_3& \cdots & h_mw_m\\
w_3 & w_2 & h_3w_3& \cdots & h_mw_m\\
\vdots & \vdots & \vdots &\ddots & \vdots\\
w_m & w_m& w_m & \cdots &h_mw_m \end{bmatrix}.$$
\item[(ii)] $\Rr^{-1}$ is a tridiagonal and irreducible matrix.
\end{itemize}

Moreover, 
$${\rm det}\,\Rr=h_1v_1w_m\prod\limits_{s=2}^mh_s(v_sw_{s-1}-v_{s-1}w_s).$$
\end{lem}

\begin{thm}
\label{Jacobi:inversa} 
The matrix $\JJ(a,b,c)$ is invertible if and only if $D_{_\JJ}\not=0$, 
and in that case, the entries of its inverse $\Rr$ are explicitly given by
$$r_{ks}= \dfrac{1}{D_{_\JJ}}\left\{ \begin{array}{cl}
\Big(\prod\limits_{j=k}^{s-1}a(j)\Big)\Phi_{_\JJ}(k)\Psi_{_\JJ}(s), & \hbox{if $0\le k\le s\le n+1$},\\[3ex]
\Big(\prod\limits_{j=s}^{k-1}c(j)\Big)\Phi_{_\JJ}(s)\Psi_{_\JJ}(k), & \hbox{if $0\le s\le k\le n+1$}. 
\end{array}\right.$$
Moreover, 
$${\rm det}\,\Rr=-D_{_\JJ}^{-1}.$$
\end{thm}

\begin{pf}
The first part is consequence of Corollary \ref{Poisson:resolvente} and taking into account the identities from Lemma \ref{phi_jacobi}. Then, for any $k,s=0,\ldots,n+1$, we obtain
$$\dfrac{\rho(s)\Big(\prod\limits_{s=0}^{n-1}c(s)\Big)}{a(0)\Big(\prod\limits_{s=1}^{\min\{k,s\}-1}a(s)\Big)\Big(\prod\limits_{s=\max\{k,s\}}^{n-1}c(s)\Big)}=\left\{ \begin{array}{cl}
\prod\limits_{j=k}^{s-1}a(j), & \hbox{si $k\le s$},\\[3ex]
\prod\limits_{j=s}^{k-1}c(j), & \hbox{si $k\ge s$}.
\end{array}\right.$$

To prove the formula for the determinant of $\Rr$, we apply Lemma \ref{Green:generalizada} with $h_j=h\rho(j)$, $h=a(0)^{-1}D_{_\JJ}^{-1}\prod\limits_{s=0}^{n-1}c(s)$, $v_j=\Phi^{a,b,c}_{\mathfrak{c}_1,\mathfrak{c}_2}(j)$ and $w_j=\Psi^{a,b,c}_{\mathfrak{c}_1,\mathfrak{c}_2}(j)$, which implies 
$$\begin{array}{rl}
h_s(v_sw_{s-1}-v_{s-1}w_s)=&\hspace{-.25cm}-h\rho_{a,c}(s)w[\Phi^{a,b,c}_{\mathfrak{c}_1,\mathfrak{c}_2},\Psi^{a,b,c}_{\mathfrak{c}_1,\mathfrak{c}_2}](s-1)\\[1ex]
=&\hspace{-.25cm}-hc(s-1)^{-1}a(s-1)\rho_{a,c}(s-1)w[\Phi^{a,b,c}_{\mathfrak{c}_1,\mathfrak{c}_2},\Psi^{a,b,c}_{\mathfrak{c}_1,\mathfrak{c}_2}](s-1)\\[1ex]
=&\hspace{-.25cm}
-hc(s-1)^{-1}a(0)w[\Phi^{a,b,c}_{\mathfrak{c}_1,\mathfrak{c}_2},\Psi^{a,b,c}_{\mathfrak{c}_1,\mathfrak{c}_2}](0)\\[1ex]
=&\hspace{-.25cm}hc(s-1)^{-1}a(0)D_{_\JJ}
\Big(\prod\limits_{s=0}^{n-1}c(s)\Big)^{-1}=c(s-1)^{-1}
,\end{array}$$
for any $s=1,\ldots,n+1$. Therefore
$${\rm det}\,\Rr=h_{0}v_{0}w_{n+1}\prod\limits_{s=1}^{n+1}c(s-1)^{-1}=-D_{_\JJ}^{-1}.$$\qed\end{pf}

Although the expression of the inverse of $\JJ(a,b,c)$ in terms of solutions of initial and final value problems is well known, see \cite{FP05,M01}, the above--explained proposal has the novelty of computing such solutions explicitly. On the other hand, the formula for the determinant of $\Rr$ appears to be new, probably because this is the first study on the inversion of matrices from an algebraic point of view, particularly based on the properties of difference equations.

We end this section particularizing the last results for a Jacobi matrix $\JJ(a,b,c)$ with constant diagonals except for the first and the last row, that is $a(j)=\alpha\not=0$, $b(j)=\beta$, $j=1,\ldots,n$,  and $c(j)=\gamma\not=0$, $j=0,\ldots,n-1$, and also for the easiest case when $\JJ(a,b,c)$ is also a Toeplitz matrix, so then has the three main diagonals completely constant. In both cases, the Schr\"odinger equation corresponds to a second order linear difference equation with constant coefficients (in the first case, the first and the last row are related to the boundary conditions), so its solution can be expressed in terms of Chebyshev polynomials, a known result that can be consulted in \cite[Theorem 2.4]{ABD05} or \cite[Theorem 2.4]{EJ16}. Of course, this result coincides with the one showed below when we use Chebyshev functions $P_k(x,y)$ valued in constant sequences $x(j)=x$ and $y(j)=y\not=0$, $j=1,\ldots,n$, so then Equation (\ref{Ch}) become Chebyshev polynomials of second kind,
$$P_{-1}(x,y)=0,\hspace{.5cm}P_0(x,y)=1\hspace{.5cm}
\hbox{and}\hspace{.5cm}P_k(x,y)=y^{\frac{k}{2}}U_k\Big(\frac{x}{2\sqrt{y}}\Big).
$$
If we consider $q=\dfrac{\beta}{2\sqrt{\alpha \gamma}}$, then
\begin{equation}
\label{base:Jacobic}
\small\begin{array}{rlr}
\Phi_{_\JJ} (0)=&\hspace{-.25cm}1, & \\[1ex]
\Phi_{_\JJ} (k)=&\hspace{-.25cm} (\sqrt{\alpha \gamma})^{k-2}\Big[b(0)\sqrt{\alpha \gamma} \,U_{k-1}(q)-a(0)\gamma U_{k-2}(q)\Big],&\hspace{-.25cm}  k=1\ldots,n+1,\\[2ex]
\Psi_{_\JJ} (k)=&\hspace{-.25cm}  (\sqrt{\alpha \gamma})^{n-k-1}\Big[b(n+1) \sqrt{\alpha \gamma}\, U_{n-k}(q)-c(n)\alpha U_{n-k-1}(q)\Big] , & \hspace{-.25cm} k=0,\ldots,n,\\[2ex]
\Psi_{_\JJ} (n+1)=&\hspace{-.25cm}1, &\end{array}
\end{equation}
and $D_{_\JJ}=d_{\JJ}(\sqrt{\alpha \gamma})^{n-2}$  where 
\begin{equation}
\label{determinant:Jacobic}
\begin{array}{rl}
d_{_\JJ}=&\hspace{-.25cm} b(0)\sqrt{\alpha \gamma}\,\Big[b(n+1)\sqrt{\alpha \gamma}\,U_{n}(q) -c(n)\alpha U_{n-1}(q)\Big]\\[1ex]
-&\hspace{-.25cm}a(0)\gamma \Big[b(n+1)\sqrt{\alpha \gamma}\,U_{n-1}(q)-c(n)\alpha  U_{n-2}(q) \Big]. 
\end{array}
\end{equation}

The next result corresponds to the first case, a Jacobi matrix with constant diagonals except for the first and the last row, and is a straightforward consequence of Theorem \ref{Jacobi:inversa} using Equations (\ref{base:Jacobic}) and (\ref{determinant:Jacobic}).

\begin{corollary}
\label{Jacobi:inversa-constante} If $a(j)=\alpha\not=0$, $b(j)=\beta$, $j=1,\ldots,n$, $c(j)=\gamma\not=0$, \linebreak$j=0,\ldots,n-1$,
then $\JJ(a,b,c)$ is invertible if and only if  $d_{_\JJ}\not=0$, 
and in that case the entries of its inverse $\Rr$ are explicitly given by
$$r_{ks}= \dfrac{1}{d_{_\JJ}(\sqrt{\alpha \gamma})^{n-2}}\left\{ \begin{array}{cl}
a(0)\alpha^{s-1}\Phi_{_\JJ}(0)\Psi_{_\JJ}(s), & \hbox{si $0=k\le s\le n+1$},\\[1ex]
\alpha^{s-k}\Phi_{_\JJ}(k)\Psi_{_\JJ}(s), & \hbox{si $1\le k\le s\le n+1$},\\[1ex]
\gamma^{k-s}\Phi_{_\JJ}(s)\Psi_{_\JJ}(k), & \hbox{si $0\le s\le k\le n$},\\[1ex]
c(n)\gamma^{n-s}\Phi_{_\JJ}(s)\Psi_{_\JJ}(n+1), & \hbox{si $0\le s\le k=n+1$}. 
\end{array}\right.$$
Moreover,
\vspace{-.2cm}$${\rm det}\,\Rr=-\dfrac{1}{d_{_\JJ}(\sqrt{\alpha\gamma})^{n-2}}.$$
\end{corollary}

Finally, the two last results showed above correspond to Jacobi and Toeplitz matrices.

\begin{corollary}
\label{Jacobi:inversa-Toeplitz} If $\alpha\gamma\not=0$, the  Jacobi and Toeplitz matrix of size $n+2$
$$\JJ(\alpha,\beta,\gamma)=\begin{bmatrix}
\beta & -\alpha & 0&\cdots & 0 & 0\\
-\gamma & \beta& -\alpha & \cdots & 0& 0\\
0 & -\gamma  & \beta& \cdots & 0 & 0\\
\vdots & \vdots &\vdots & \ddots & \vdots & \vdots \\
0 & 0 & 0 & \cdots & \beta & -\alpha\\
0 & 0& 0& \cdots  & -\gamma  & \beta\end{bmatrix}$$
is invertible if and only if 
$$\beta\not=2\sqrt{\alpha\gamma}\cos\left(\frac{k\pi}{n+3}\right),\hspace{.5cm}k=1,\ldots,n+2,$$
and then, the entries of the inverse of $\JJ(\alpha,\beta,\gamma)$ are explicitly given by
$$r_{ks}= \dfrac{1}{U_{n+2}(q)}\left\{ \begin{array}{cl}
\alpha^{s-k}(\sqrt{\alpha \gamma})^{k-s-1} U_{k}(q) U_{n-s+1}(q), & \hbox{if $0\le k\le s\le n+1$},\\[2ex]
\gamma^{k-s}(\sqrt{\alpha \gamma})^{s-k-1} U_{s}(q)U_{n-k+1}(q), & \hbox{if $0\le s\le k\le n+1$},
\end{array}\right.$$
where $q=\dfrac{\beta}{2\sqrt{\alpha \gamma}}$.

\medskip

Moreover,
$${\rm det}\,\Rr=\dfrac{1}{\sqrt{\alpha\gamma})^{(n+2)}U_{n+2}(q)}.$$
\end{corollary}

\begin{pf}  
All the results are consequence of Theorem \ref{Jacobi:inversa} by imposing in Equations (\ref{base:Jacobic}) and (\ref{determinant:Jacobic}) the identities
$$a(0)=-\alpha,\hspace{.15cm}b(0)=b(n+1)=\beta,\hspace{.15cm}c(n)=-\gamma.$$
Then, 
$$
d_{_\JJ}=-\alpha^2\gamma^{2}U_{n+2}(q),
$$ 
so $d_{_\JJ}\neq 0$ if and only if $q$ is not a zero of the polymonial $U_{n+2}(x)$; that is, if and only if $q\not=\cos\big(\frac{k\pi}{n+3}\big)$, $k=1,\ldots,n+2$, see \cite{MH03}. Moreover, the expression for the determinant follows. 

On the other hand,
$$
\Phi_{_\JJ}(k)=(\sqrt{\alpha \gamma})^{k} U_{k}(q), \hspace{.25cm}
\Psi_{_\JJ}(k)= -(\sqrt{\alpha \gamma})^{n-k+1} U_{n-k+1}(q) 
$$
for any $k=0,\ldots,n+1$, that leads to the given expressions for the inverse entries. \qed 
\end{pf}

A more detailed proof of the above result for Jacobi and Toeplitz matrices can be consulted in \cite{EnJi17}. Besides, the expression obtained for the matrix inverse of this kind of matrices coincides with that published by Fonseca and Petronilho in \cite[Corollary 4.1]{FP01} and \cite[Equation 4.26]{FP05}. 

\begin{corollary}
\label{Jacobi:inversa-Toeplitz:sym} If $\alpha\not=0$, the symmetric Jacobi and Toeplitz matrix of order $n+2$
$$\JJ(\alpha,\beta)=\begin{bmatrix}
\beta & -\alpha & 0&\cdots & 0 & 0\\
-\alpha & \beta& -\alpha & \cdots & 0& 0\\
0 & -\alpha  & \beta& \cdots & 0 & 0\\
\vdots & \vdots &\vdots & \ddots & \vdots & \vdots \\
0 & 0 & 0 & \cdots & \beta & -\alpha\\
0 & 0& 0& \cdots  & -\alpha  & \beta\end{bmatrix}$$
is invertible if and only if
$$\beta\not=2\alpha\cos\left(\frac{k\pi}{n+3}\right),\hspace{.5cm}k=1,\ldots,n+2,$$
and then, the entries of the inverse of $\JJ(\alpha,\beta)$ are explicitly given by
$$r_{ks}= \dfrac{U_{\min\{k,s\}}\big(\frac{\beta}{2\alpha}\big) U_{n-\max\{k,s\}+1}\big(\frac{\beta}{2\alpha}\big)}{\alpha U_{n+2}\big(\frac{\beta}{2\alpha}\big)} , \hspace{.25cm}k,s=0,\ldots,n+1.$$
Moreover,
$${\rm det}\,\Rr=\dfrac{1}{\alpha^{(n+2)}U_{n+2}\big(\frac{\beta}{2\alpha}\big)}.$$
\end{corollary}
\medskip
The expression for the  inverse of a symmetric Jacobi and Toeplitz matrix is well--known, see for instance \cite[Corollary 4.2]{FP01} and the references of this article.

\bigskip
This work has been partly supported by the Spanish Program I+D+i (Mi\-nis\-te\-rio de Econom\'ia y Competitividad) under projects MTM2014-60450-R and MTM2017-85996-R.

\section*{References}

\bibliography{mybibfile}

\end{document}